\documentclass[12pt]{article}
\usepackage{epsfig}
\usepackage{graphicx}

\usepackage{amsfonts}
\usepackage{amssymb}
\usepackage{amsthm}

\def\RR{\leavevmode\hbox{$\rm I\!R$}}

\newtheorem{propo}{Proposition}[section]
\newtheorem{lemma}[propo]{Lemma}

\def\pmb#1{\setbox0=\hbox{$#1$}%
             \kern-.033em\copy0\kern-\wd0
             \kern+.011em\copy0\kern-\wd0
             \kern+.011em\copy0\kern-\wd0
             \kern+.011em\copy0\kern-\wd0
             \kern+.011em\copy0\kern-\wd0
             \kern+.011em\copy0\kern-\wd0
             \kern+.011em\copy0\kern-\wd0
             \kern-.055em\raise+.015em\copy0\kern-\wd0
             \kern+.011em\raise+.015em\copy0\kern-\wd0
             \kern+.011em\raise+.015em\copy0\kern-\wd0
             \kern+.011em\raise-.015em\copy0\kern-\wd0
             \kern+.011em\raise-.015em\copy0\kern-\wd0
             \kern-.022em\copy0\kern-\wd0\raise-.015em\box0}

\def\BBox{\kern  -0.2cm\hbox{\vrule width 0.15cm height 0.3cm}}

\def\.#1{{\buildrel .\over {#1}}}

\def\proof{{\noindent\bf Proof.}\hskip 0.3truecm}
\def\BBox{\kern  -0.2cm\hbox{\vrule width 0.15cm height 0.3cm}}

\begin{document}

\title{Delaunay Surfaces}

\author{Enrique Bendito\footnote{\small e-mail: enrique.bendito@upc.edu}, Mark J. Bowick{\small **} and Agust\'\i n Medina{\small*}\\
{*\em Departament Matem\`atica Aplicada III}\\
{Universitat Polit\`ecnica de Catalunya,}\\
{Barcelona, Spain.}\\
{**\em Physics Department}\\
{Syracuse University}\\
{Syracuse, NY 13244-1130, USA}}

\date{}

\maketitle

\begin{abstract}

We derive parametrizations of the Delaunay constant mean curvature surfaces of revolution that follow directly from 
parametrizations of the conics that generate these surfaces via the corresponding roulette. This uniform treatment 
exploits the natural geometry of the conic (parabolic, elliptic or hyperbolic) and leads to simple expressions for 
the mean and Gaussian curvatures of the surfaces as well as the construction of new surfaces.
\end{abstract}
\section{Preliminaries}

The surfaces of revolution with constant mean curvature (CMC) were introduced and completely characterized by C. Delaunay more than a 
century ago \cite{D41}. Delaunay's formulation of the problem leads to a non--linear ordinary differential equation  involving the 
radius of curvature of the plane curve that generates the surface, which can be also characterized variationally as the surface of revolution 
having a minimal lateral area with a fixed volume (see \cite{E87}). Delaunay showed that the above differential equation arises 
geometrically by rolling a conic along a straight line without slippage. The curve described by a focus of the conic, 
the \emph{roulette of the conic}, is then the meridian of a surface of revolution with constant mean curvature, where the 
straight line is the axis of revolution. These CMC surfaces of revolution are called \emph{Delaunay surfaces}.
Apart from the elementary cases of spheres and cylinders, there are three classes of Delaunay surfaces,  the \emph{catenoids}, 
the \emph{unduloids} and the \emph{nodoids}, corresponding to the choice of conic as a parabola, an ellipse or a hyperbola, respectively.

Traditionally the roulettes have been characterized using polar coordinates centered at the focus of the 
conic~\cite{B90,K80,HY81,E87,M02,R05,HMO07,KP08}. The methods emplyed in these papers are based on solving 
certain ordinary differential equations that, in one way or another, depend on the variational characterization of the CMC surfaces. 
Although Ref.~\cite{B90} does suggest the possibility of using the cartesian coordinates of the roulettes with the tangent to the 
conic as the abscissa, this idea is never developed. 
  
Here we obtain parametrizations of the roulettes, and therefore of the corresponding Delaunay surfaces, directly from the parametrizations 
of the conics. This leads directly to concise expressions for all the key differential geometric characteristics of Delaunay surfaces.
In our approach the unduloid is described with trigonometric functions, whereas the catenoid and the nodoid are described with 
hyperbolic functions. This yields simple expressions for the Gaussian curvature, total curvature and mean 
curvature as well as the length of roulettes. The mean curvature of an unduloid, in particular, is given by the inverse of the distance 
between the vertices of the corresponding ellipse, whereas the mean curvature of a nodoid is given by minus the inverse of the distance 
between the vertices of the corresponding hyperbola. The parametrizations presented here also give rise to a straightforward construction 
of nodoids, both when viewed as simple parts (generated by a focus) or when they are composed of several individual parts or a periodic 
repetition of simple parts.

For the sake of  completeness, we finish this section by  presenting some well-known results about regular surfaces of revolution, 
as well as a very simple proof of the Gauss-Bonett theorem for this class of surfaces.

Let $f,g:[t_1,t_2]\longrightarrow\mathbb{R}$ be smooth functions with $f>0$, and $S$ the  surface of revolution parametrized by 
$\pmb x:[t_1,t_2]\times[v_1,v_2]\longrightarrow\mathbb{R}^3$ %
$$\pmb x(t,v)=\Big(f(t)\cos(v),f(t)\sin(v),g(t)\Big).$$
The coefficients of the first and second fundamental forms of $S$ are given by
$$\begin{array}{rlll}
E=&\hspace{-.25cm}\langle {\pmb x}_t,{\pmb x}_t\rangle=(f')^2+(g')^2,\hspace{.2cm} & F=
\langle {\pmb x}_t,{\pmb x}_v\rangle= 0, \hspace{.2cm} &
 G=\langle {\pmb x}_v,{\pmb x}_v\rangle=f^2;\\[2ex]
L=&\hspace{-.25cm}\langle {\pmb x}_{tt},{\pmb n}\rangle=\displaystyle\frac{f'g''-f''g'}{\Big((f')^2+(g')^2\Big)^{1/2}},\hspace{.2cm} & 
M=\langle {\pmb x}_{tv},{\pmb n}\rangle=0,\hspace{.4cm} & N=
\langle {\pmb x}_{vv},{\pmb n}\rangle=\displaystyle\frac{fg'}{\Big((f')^2+(g')^2\Big)^{1/2}},
\end{array}$$
where
$${\pmb n}=\displaystyle\frac{{\pmb x}_t\times {\pmb x}_v}{|{\pmb x}_t\times {\pmb x}_v|}=\displaystyle\frac{1}{\Big((f')^2+(g')^2\Big)^{1/2}}
\Big(- g'\cos(v),-g'\sin(v),f'\Big),$$
is the unit normal to $S$. The Gaussian curvature is given by
$$K=\displaystyle\frac{LN}{EG}=\frac{g'\big(f'g''-f''g'\big)}{f\Big((f')^2+(g')^2\Big)^2},$$
whereas the  mean curvature, $H$, is given by
$$2H=k_1+k_2=\displaystyle\frac{L}{E}+\frac{N}{G}=\frac{f'g''-f''g'}{\Big((f')^2+(g')^2\Big)^{3/2}}+
\frac{g'}{f\Big((f')^2+(g')^2\Big)^{1/2}},$$
where $k_1$ and $k_2$ are the two principal curvatures.

Now consider a curve $\alpha$ in $S$ parametrized by the arc length $s$. For any point $p=\alpha(s)$ we choose a vector 
$\pmb u(s)$ in the tangent space at $p$ such that $\{\alpha'(s),\pmb u(s)\}$ is a positively oriented orthonormal basis 
of the tangent space at $p$; that is, $\alpha'(s)\times \pmb u(s)=\pmb n(\alpha(s))$. The geodesic curvature $k_g(s)$ of 
$\alpha$ at $s$ is then given by
$$k_g(s)=\langle \displaystyle\alpha''(s), \pmb u(s)\rangle.$$
If the curve is chosen to be a meridian of $S$, $\alpha(t)=\pmb x(t,v_0)$, then $k_g=0$, whereas if it is a parallel, 
$\beta(v)=\pmb x(t_0,v)$, then its geodesic curvature is given by
$$k_g(t_0)=\displaystyle\frac{f'(t_0)}{f(t_0)\Big((f'(t_0))^2+(g'(t_0))^2\Big)^{1/2}}.$$
\begin{lemma} If  $C_1$ and $C_2$ are the boundary parallels of $S$ with the orientation induced by $S$ then
$$\int_SKd\sigma+\int_{C_1}k_g(t_1)d\ell+\int_{C_2}k_g(t_2)d\ell=0.$$
\end{lemma}
\proof Observe first that
  $\displaystyle k_g(t)|\pmb x_v|=\frac{f'(t)}{\Big((f'(t))^2+(g'(t))^2\Big)^{1/2}}$ and hence
$$\Big(k_g |\pmb x_v|\Big)'= \frac{g' \Big(f'' g' -f' g'' \Big)}{\Big((f' )^2+(g' )^2\Big)^{3/2}}=-K |\pmb x_t|\,|\pmb x_v|.$$

On the other hand,
$$\begin{array}{rl}
\displaystyle\int_SKd\sigma=&\hspace{-.25cm}\displaystyle\int_{v_1}^{v_2}\int_{t_1}^{t_2}K|\pmb x_t|\,|\pmb x_v|dtdv= 
-\int_{v_1}^{v_2}\int_{t_1}^{t_2}\Big(k_g |\pmb x_v|\Big)'dtdv\\[3ex]
=&\hspace{-.25cm}\displaystyle \int_{v_1}^{v_2} k_g (t_1)|\pmb x_v(t_1,v)|dv-\int_{v_1}^{v_2} k_g (t_2)|\pmb x_v(t_2,v)|  
dv=-\int_{C_1}k_g(t_1)d\ell-\int_{C_2}k_g(t_2)d\ell.\end{array}$$\qed

If $v_2-v_1=2\pi$, then $S$ is homeomorphic to an annulus with  null Euler characteristic. Otherwise, it is  a simple region whose 
Euler characteristic equals $1$ and the sum  of the angles at the four boundary vertices, formed by the tangents to the boundary 
curves oriented with the orientation induced by $S$, equals $2\pi$.

\section{The Roulettes of the conics} 
\begin{figure}[htb]
\begin{center}
\scalebox{0.39}{{\includegraphics{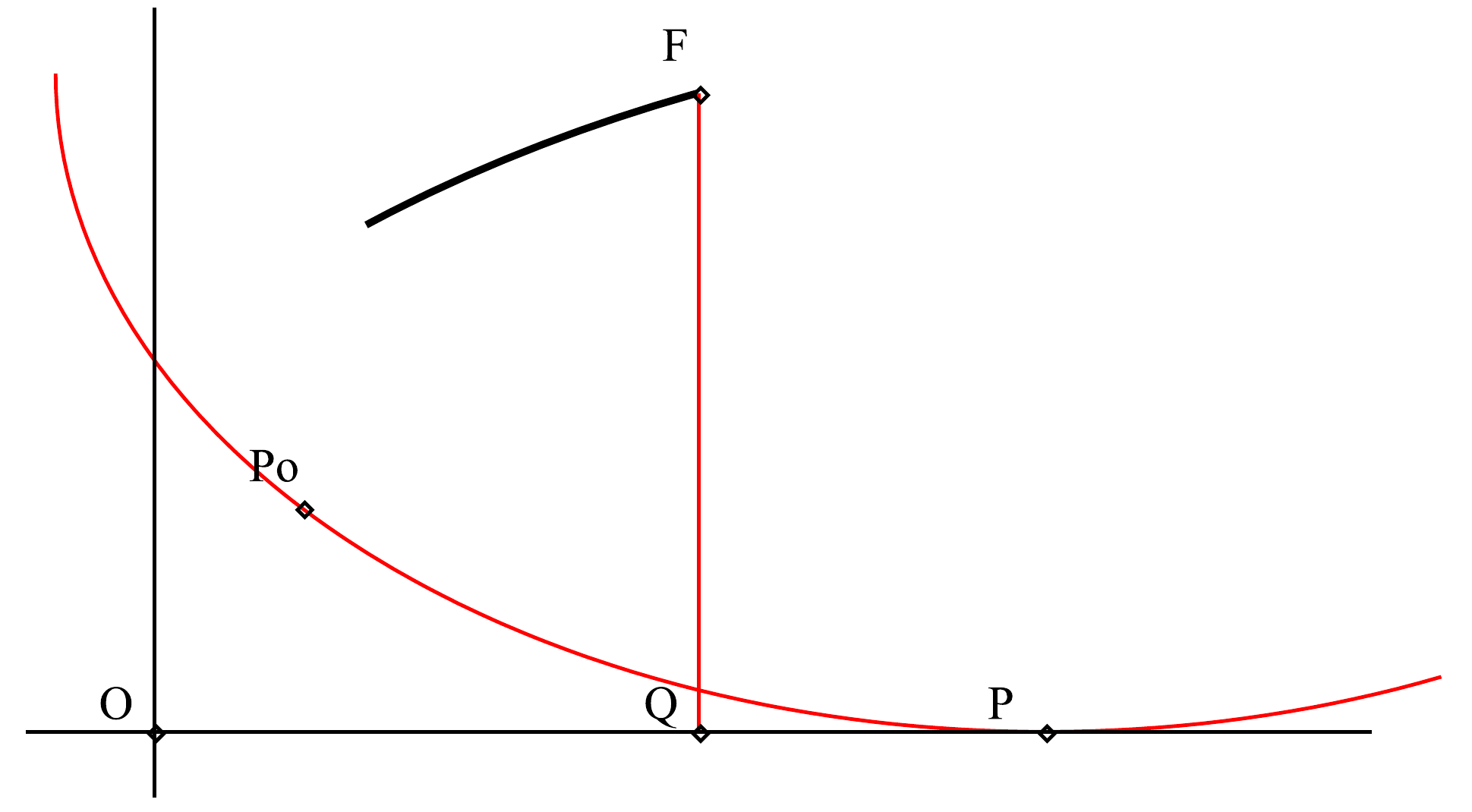}}}\hspace{1.5cm}\scalebox{0.32}
{{\includegraphics{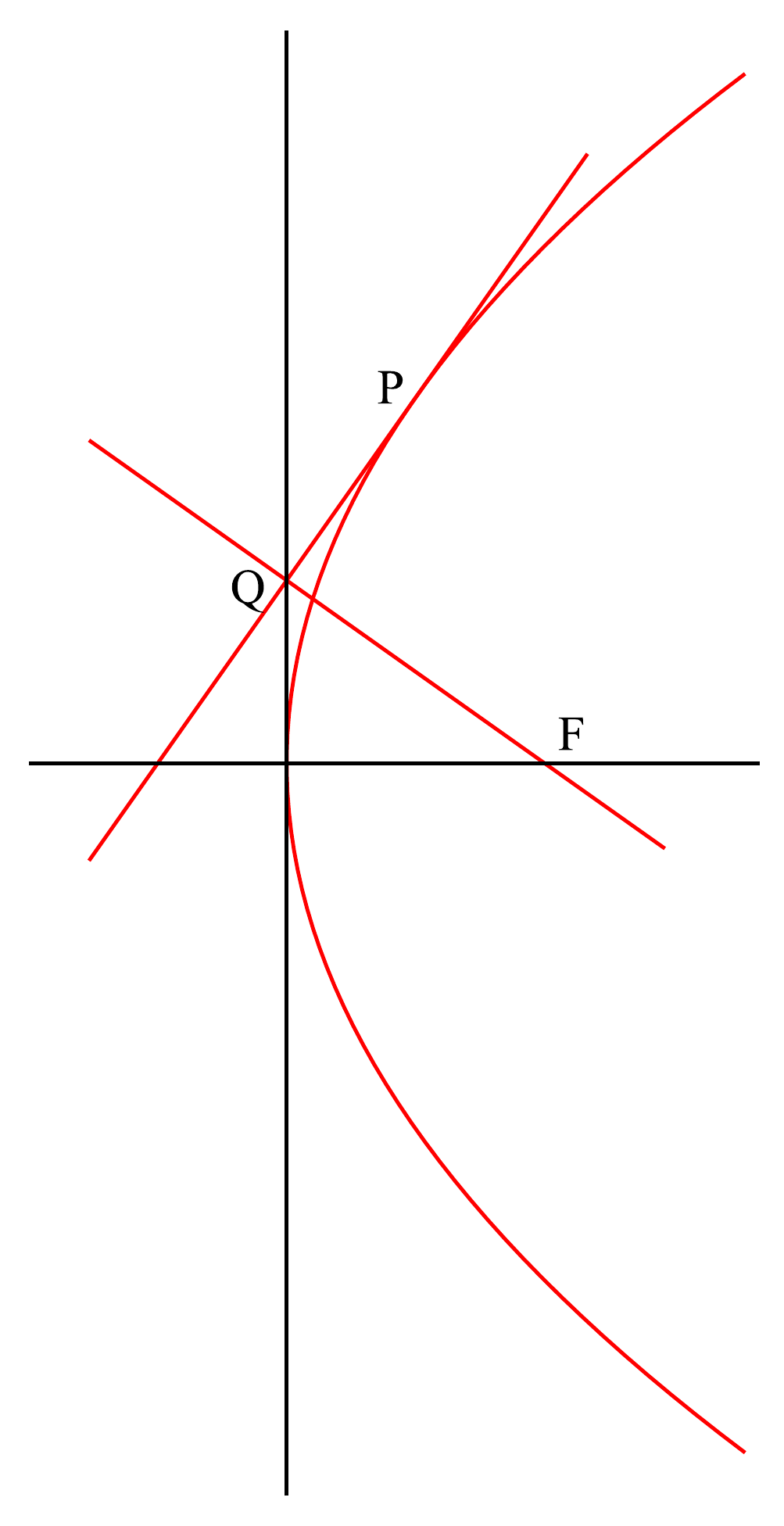}}}
\end{center}
\caption{\label{parabola} Roulette (left) and Parabola (right).}
\end{figure}

When a curve rolls, without slipping, on a fixed curve, each  point of the rolling 
curve traces another curve known as a {\it roulette}. In Figure \ref{parabola} (left), we show the trace generated by the 
point $F=(F_1,F_2)$, associated with a given curve $C$, when it rolls on a straight line. The abscissa $F_1$ of the roulette 
coincides with $Q_1$, i.e., the length of the arc of the curve from $P_o$ to $P$, minus the value $P_1-Q_1$. We can interchange
the roles of the conic and its tangent by considering a fixed conic and a moving tangent. The locus of the points $Q$ thus 
obtained is called the {\it pedal curve} of $C$ with respect to $F$. Here we are interested in the roulettes generated by the conic 
foci when they roll over a tangent line.

The parametric description of the {\it parabola}, the {\it ellipse} and the {\it hyperbola} are given respectively by
$$\begin{array}{rl}
\displaystyle
 \alpha(t)=&\hspace{-.2cm}\left(b\sinh^2(t),2b\sinh(t)\right),\\[2ex]
\displaystyle
 \beta(t)=&\hspace{-.2cm}\Big(a\cos(t),b\sin(t)\Big),\\[2ex]
\displaystyle
 \gamma(t)=&\hspace{-.2cm}\Big(a\cosh(t),b\sinh(t)\Big),
\end{array}$$
where $a, b>0$ and $t\in [t_1,t_2]\subset\mathbb{R}$.

In Figure \ref{parabola} (right), we show a parabola, with the perpendicular from $F$ onto the tangent of the parabola at $P$. 
This situation is equivalent to the general case and we apply the same criteria to give a parameterization of this roulette.  

The arc length for the parabola from  $t_0$ is
$$s=\int_{t_0}^t|\alpha '(u)|du=b(t+\sinh(t)\cosh(t)).$$
The length of the segment $\overline{PQ}$ is $b\sinh(t)\cosh(t)$ yielding an abscissa 
$$
  g_c(t)=s-b\sinh(t)\cosh(t)=bt.
$$
Computing $f_c(t)$, the length of the segment $\overline{FQ}$, shows that the roulette associated with the focus of the parabola
is the {\it catenary}
$$
 A(t)=(g_c(t),f_c(t))=\left(b t,b\cosh\left(t\right)\right).
$$
Note also that
$$|A'(t)|=b\cosh(t),$$
and the arc length is given by
$$\ell_c(t)=\int_{t_0}^t|A'(z)|\,dz= b\sinh(z)\Big|_{t_0}^t.$$

\begin{figure}[htb]
\begin{center}
\scalebox{0.35}{{\includegraphics{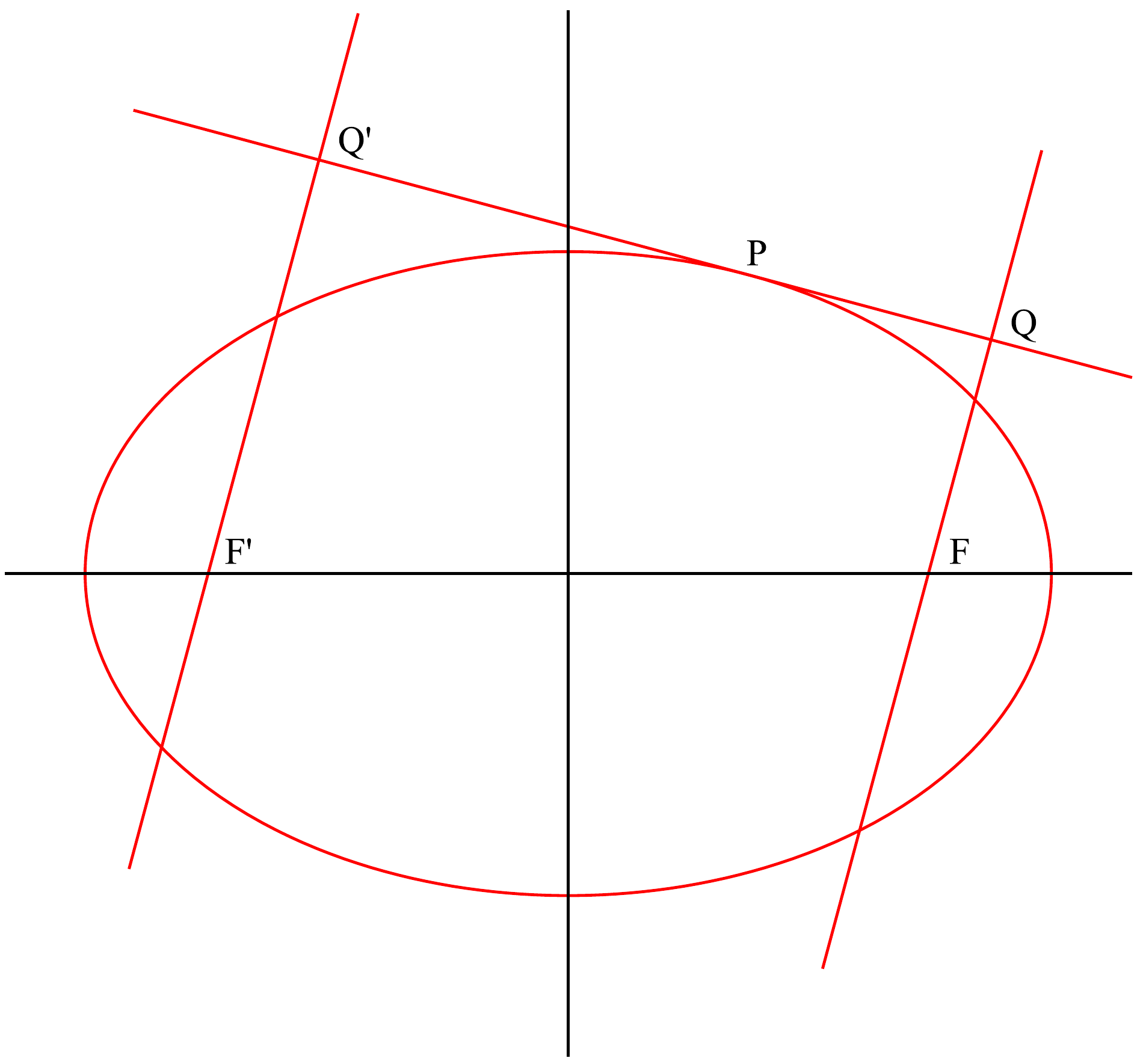}}}\hspace{1.5cm}\scalebox{0.35}
{{\includegraphics{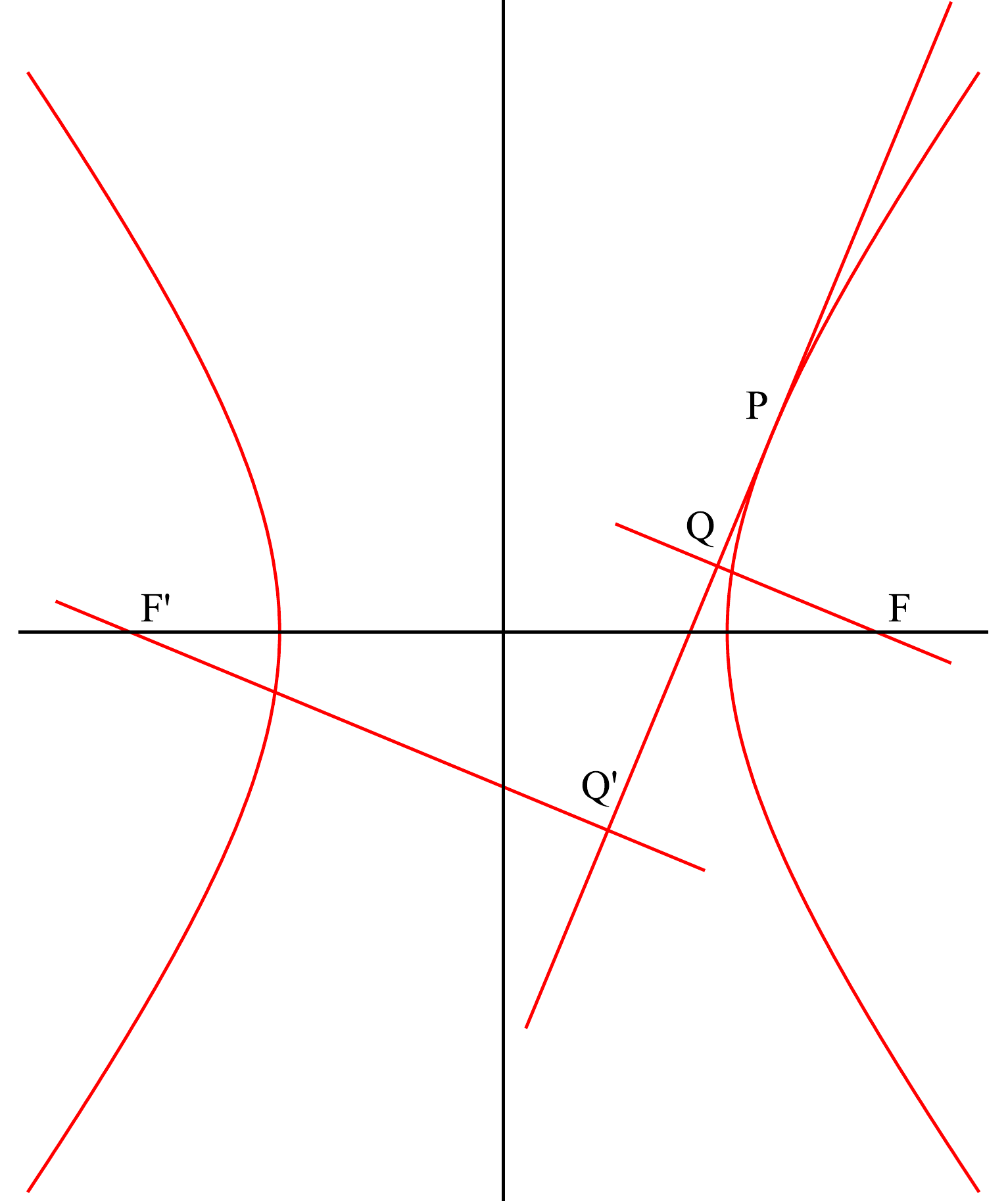}}}
\end{center}
\caption{\label{elipse} Ellipse (left) and Hyperbola (right).}
\end{figure}

For the ellipse, Figure \ref{elipse} (left), take $b<a$ and $c=\sqrt{a^2-b^2}$. The arc length from $t_0$ to $t$ is
$$s=\int_{t_0}^t|\beta '(z)|dz=\int_{t_0}^t\sqrt{a^2-c^2\cos^2(z)}\,dz.$$
In this case two curves are generated. The first one corresponds to choosing the focus $F$ closest to the tangent. 
By computing the length of the segment $\overline{PQ}$ and substracting it from the above arc, we find that 
the abscissa is 
$$g^1_u(t)=\int_{t_0}^t\sqrt{a^2-c^2\cos^2(z)}dz-\frac{c\sin(t)\left(a-c\cos(t)\right)}
{\sqrt{a^2-c^2\cos^2(t)}}.$$
In addition, the ordinate is given by the length of the segment $\overline{FQ}$; namely
$$  f^1_u(t)=\frac{b\left(a-c\cos(t)\right)}{\sqrt{a^2-c^2\cos^2(t)}}.$$
$B_1(t)=(g^1_u(t),f^1_u(t))$ is therefore the parametrization of the roulette generated by the focus of the ellipse. 
One finds
$$|B_1'(t)|=\displaystyle\frac{ab}{a+c\cos(t)},$$
and the arc length is given by
$${\ell^1_u}(t)=\int_{t_0}^t|B_1'(z)|dz=
\left.2a\arctan\left(\sqrt{\frac{a-c}{a+c}}\tan\left(\frac{z}{2}\right)\right)\right|_{t_0}^t.$$

In the same way, if we chooose the other focus $F'$, it follows after computing the length of
$\overline{PQ'}$ that the abscissa is
$$g^2_u(t)=\int_{t_0}^t\sqrt{a^2-c^2\cos^2(z)}dz-\frac{c\sin(t)\left(a+c\cos(t)\right)}
{\sqrt{a^2-c^2\cos^2(t)}},$$
and the ordinate is the length of the segment $\overline{F'Q'}$; namely

$$  f^2_u(t)=\frac{b\left(a+c\cos(t)\right)}{\sqrt{a^2-c^2\cos^2(t)}}.$$
$B_2(t)=(g^2_u(t),f^2_u(t))$ is therefore the parametrization of the roulette generated by the focus $F'$.
One now finds 
$$|B_2'(t)|=\displaystyle\frac{ab}{a-c\cos(t)},$$
and an arc length 
$${\ell^2_u}(t)=\int_{t_0}^t|B_2'(z)|dz=
\left.2a\arctan\left(\sqrt{\frac{a+c}{a-c}}\tan\left(\frac{z}{2}\right)\right)\right|_{t_0}^t.$$
Observe in particular that 
$\arctan\left(\sqrt{\frac{a+c}{a-c}}\,\right)+ \arctan\left(\sqrt{\frac{a-c}{a+c}}\,\right)=\frac{\pi}{2}$
and then the sum of the length of the two curves for $t\in (\frac{-\pi}{2},\frac{\pi}{2})$ is $2\pi a$. 

The roulette of the focus of an ellipse is called an {\it undulary}. It is clear that we do not need to consider both foci for an ellipse.
Specifically, if we consider the ellipse described by taking $t\in [-\pi,\pi]$, the curve generated by the focus $F$ is the same as the curve 
that results from joining the two curves generated by both foci $F$ and $F'$ with $t\in[-\pi/2,\pi/2]$ only. We consider both foci to 
make manifest the constructive paralelism between these roulettes and the roulettes of the hyperbola.

Consider now the case of the hyperbola, as shown Figure \ref{elipse} (right). Taking $c=\sqrt{a^2+b^2}$, the arc length from $t_0$ to $t$ is
 $$s=\displaystyle\int_{t_0}^t|\gamma'(z)|dz=
\displaystyle\int_{t_0}^t\sqrt{c^2\cosh^2(z)-a^2}\,dz.$$
For the first roulette we consider the focus $F$ closest to the tangent. By computing the length of the segment $\overline{PQ}$
it then follows that the abscissa is
$$ g^1_n(t)=\int_{t_0}^t\sqrt{c^2\cosh^2(z)-a^2}\,dz-
\frac{c\sinh(t)\left(c\cosh(t)-a\right)}{\sqrt{c^2\cosh^2(t)-a^2}},
$$
whereas its ordinate is given by the length of $\overline{FQ}$, namely,
$$f^1_n(t)=\frac{b\left(c\cosh(t)-a\right)}{\sqrt{c^2\cosh^2(t)-a^2}}.
$$
$C_1(t)=(g^1_n(t),f^1_n(t))$ is therefore the parametrization of the roulette generated by the focus $F$. 
One finds
$$|C_1'(t)|=\displaystyle\frac{ab}{c\cosh(t)+a},$$
with arc length given by
$$\ell_1(t)=\int_{t_0}^t|C_1'(z)|dz=
\left.2a\arctan\left(\sqrt{\frac{c-a}{c+a}}\tanh\left(\frac{z}{2}\right)\right)\right]_{t_0}^t.$$
In particular, the length of $C_1$ with \,$t\in (-\infty,\infty)$ is $4a\arctan\left(\sqrt{\frac{c-a}{c+a}}\right).$
\begin{figure}[htb]
\begin{center}
\scalebox{0.35}{{\includegraphics{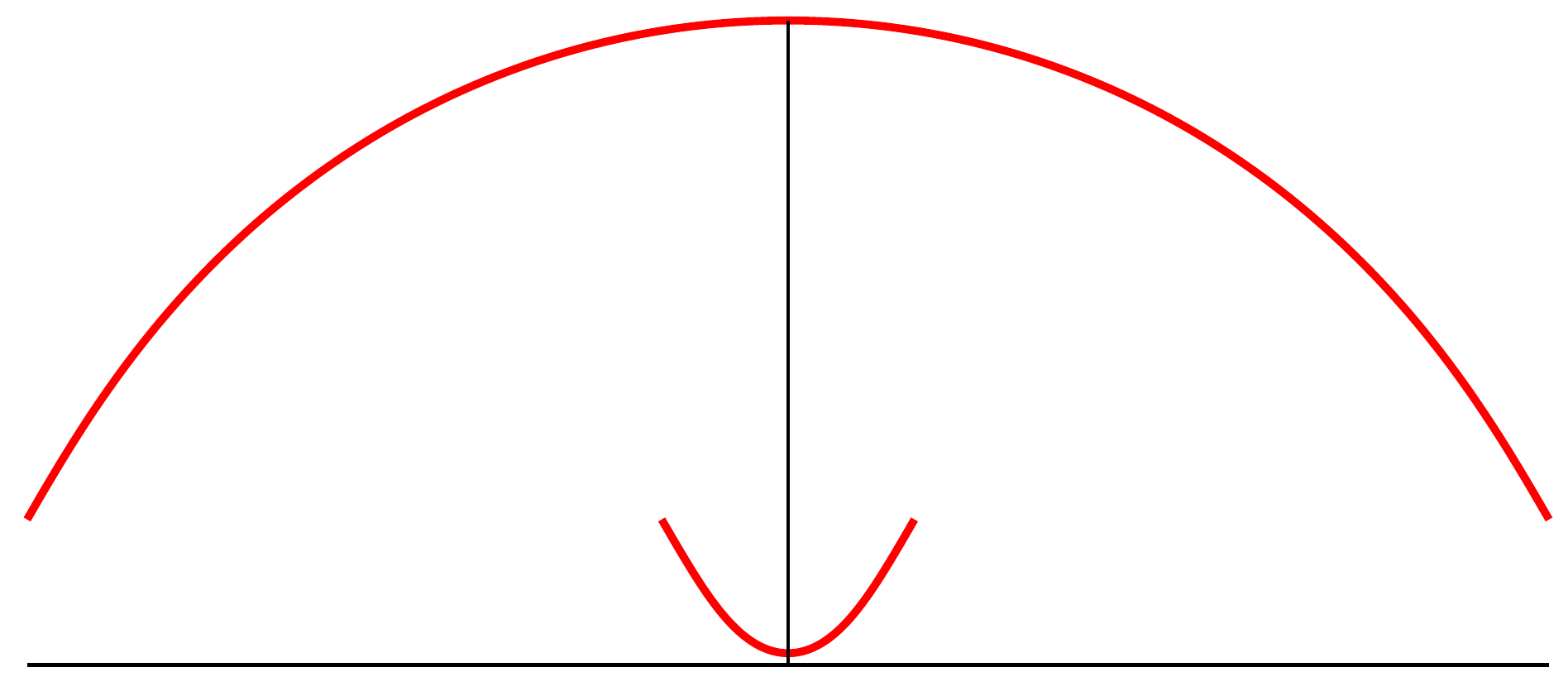}}}\hspace{1.5cm}\scalebox{0.35}
{{\includegraphics{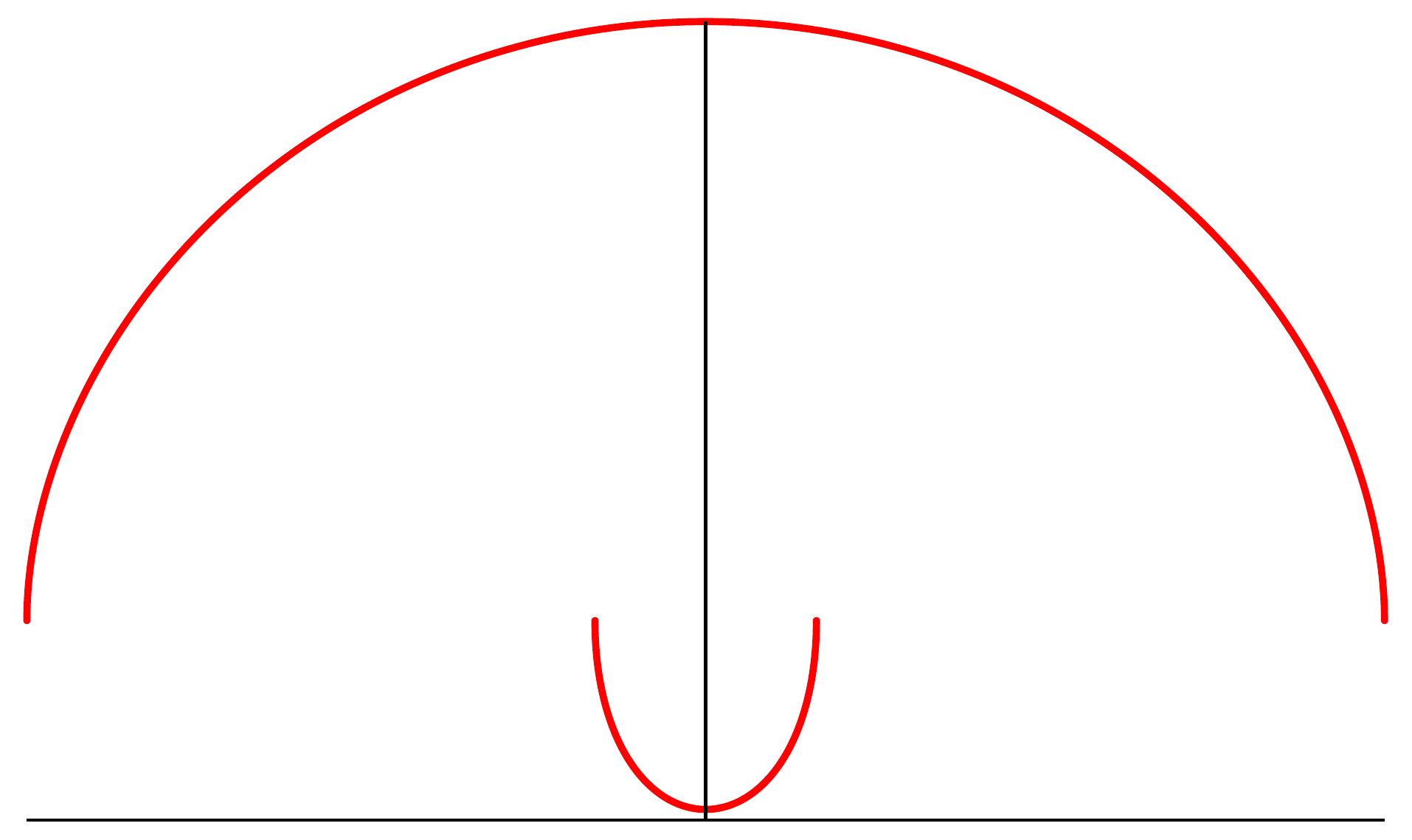}}}
\end{center}
\caption{\label{roulettes} Roulettes $B_1$ and $B_2$ (left) and roulettes $C_1$ and $C_2$ (right).}
\end{figure}

Taking the focus $F'$ instead one finds, after computing the length of $\overline{PQ'}$, that the abscissa is
$$  g^2_n(t)=\int_{t_0}^t\sqrt{c^2\cosh^2(z)-a^2}\,dz-\frac{c\sinh(t)
\left(c\cosh(t)+a\right)}{\sqrt{c^2\cosh^2(t)-a^2}},
$$
and the ordinate is the length of the segment $\overline{F'Q'}$, namely,
$$
  f^2_n(t)=\frac{b\left(c\cosh(t)+a\right)}{\sqrt{c^2\cosh^2(t)-a^2}}.
$$
$C_2(t)=(g^2_n(t),f^2_n(t))$ is therefore the parametrization of the roulette generated by the focus $F'$. 
Now 
$$|C_2'(t)|=\displaystyle\frac{ab}{c\cosh(t)-a},$$
and the arc length is given by
$${\ell_2}(t)=\int_{t_0}^t|C_2'(z)|dz=
\left.2a\arctan\left(\sqrt{\frac{c+a}{c-a}}\tanh\left(\frac{z}{2}\right)\right)\right]_{t_0}^t.$$
The sum of the length of the two branches of the curves for $t\in (-\infty, \infty)$ is again $2\pi a$. 
The roulette of the focus of a hyperbola is called a {\it nodary}.

In Figure \ref{roulettes} we display the roulettes generated by the foci of the ellipse (left) and the hyperbola (right). 
$B_1$ and $C_1$ are the curves with increasing slope, while $B_2$ and $C_2$ are the curves with decreasing slope.

\begin{figure}[htb]
\begin{center}
\scalebox{0.28}{{\includegraphics{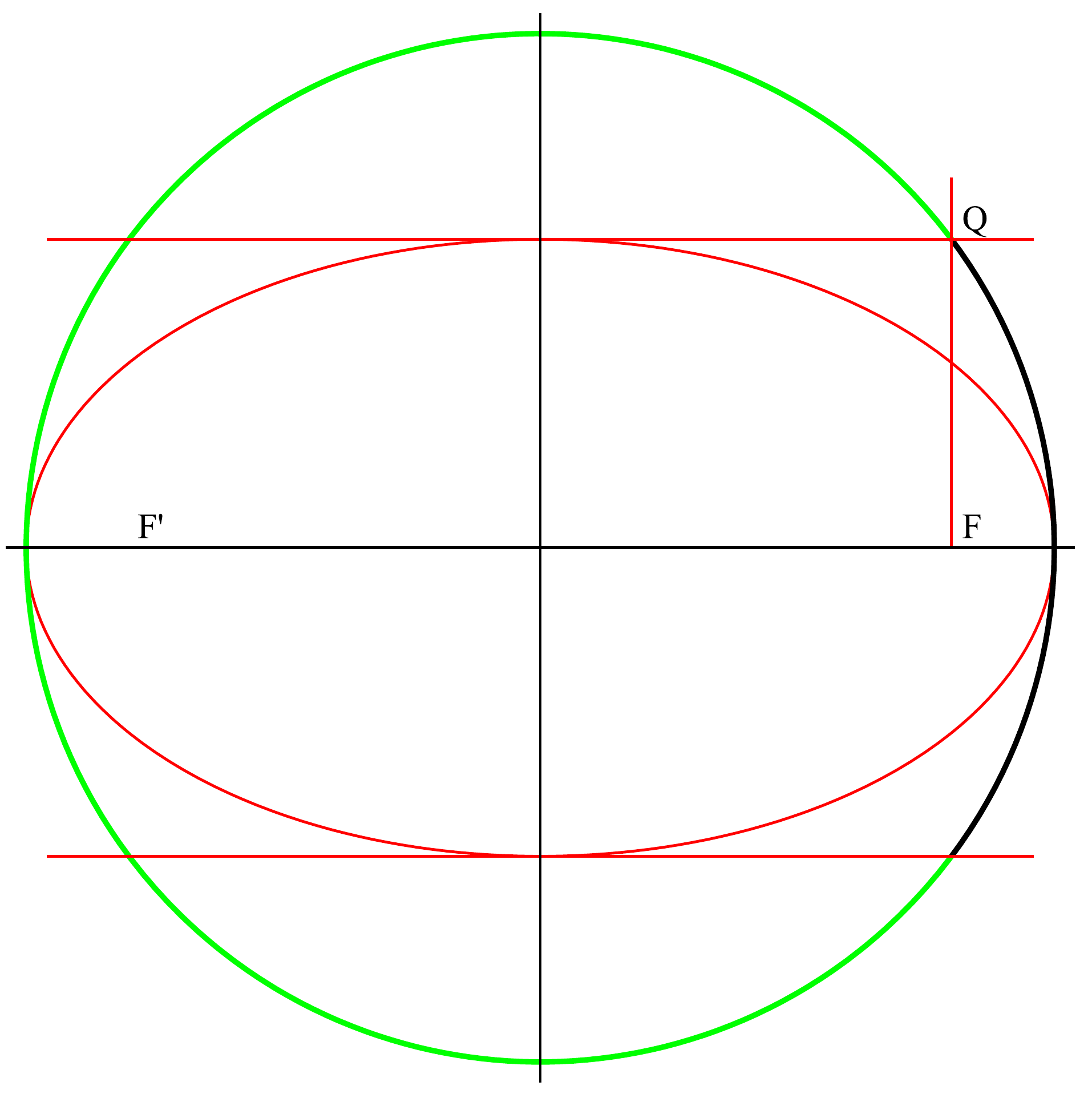}}}\hspace{1.5cm}\scalebox{0.3}
{{\includegraphics{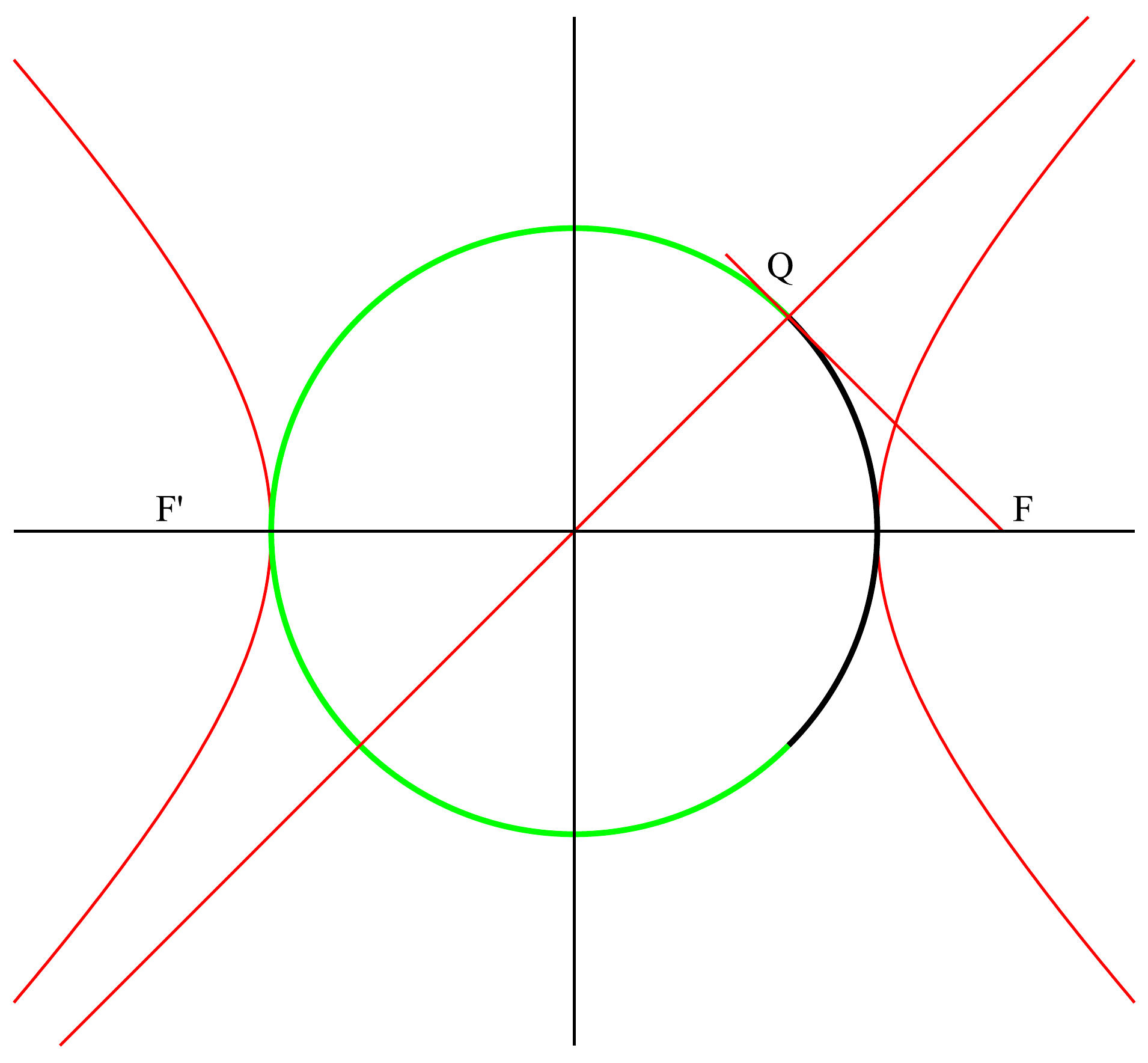}}}
\end{center}
\caption{\label{arco}  Pedal curves of the ellipse and the hyperbola. Left: arc of the circle (black), whose length coincides 
with the length of the curve $B_1$ and arc of the circle (green), whose length coincides with the length of the curve $B_2$. 
Right: arc of the circle (black), whose length coincides with the length of the curve $C_1$ and arc of the circle (green) 
whose length coincides with the length of the curve $C_2$.}
\end{figure}
In Figure \ref{arco} we display the pedal curves of the ellipse (left) and of the hyperbola (right). 
Both are circles of diameter equal to the distance between the vertices of the conic.
It is known that when a curve rolls on a straight line the arc of the roulette is equal to the corresponding arc of the pedal. 
Here we can see this directly. If we consider, for instance, the curve $C_1$ and denote $x=\arctan\left(\sqrt{\frac{c-a}{c+a}}\,\right),$
then $$\sin(2x)=\frac{b}{c} \hspace{2cm}\cos(2x)=\frac{a}{c},$$ and so 
$$ 4a \arctan\left(\sqrt{\frac{c-a}{c+a}}\,\right)= 2a\arctan\left(\frac{b}{a}\right).$$
That is, the length of the curve $C_1$ coincides with the length of the associated pedal curve {--} see Figure \ref{arco} (right in black).

\begin{figure}[htb]
\begin{center}
\scalebox{0.6}{{\includegraphics{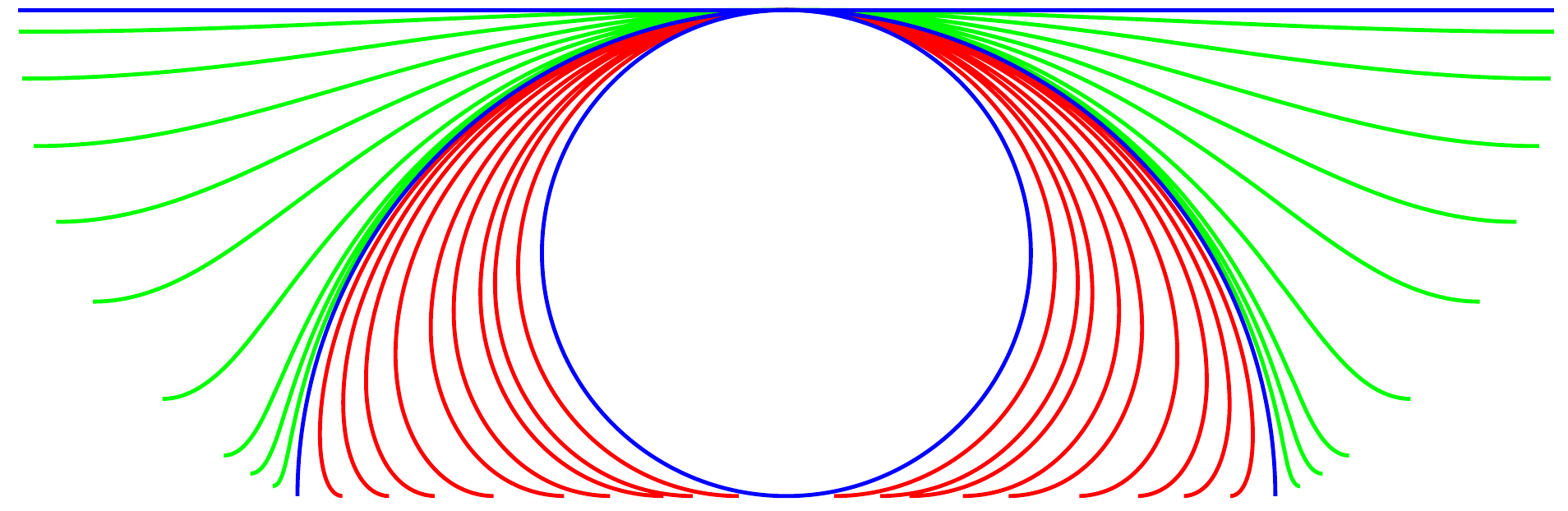}}}
\end{center}
\caption{\label{familia1} Family of roulettes $B_1\cup B_2$ and $C_1 \cup C_2$ of length $2\pi a$.}
\end{figure}

In Figure \ref{familia1} we display a family of roulettes with length $2\pi a$. Each of these rouletes was generated by a conic with 
major axis $a$. The straight segment corresponds to a curve of type $B$ in the limiting case $b=a$. This is followed by a set of elliptical 
roulettes from $b=a$ until the next limiting case $b=0$, corresponding to a semicircumference of radius $2a$. The remaining roulettes 
correspond to curves of type $C$, ranging from the limiting case $b=0$, again the same semicircumference, until $b=\infty$, which is 
the circle of radius $a$ shown in the figure.

\begin{figure}[htb]
\begin{center}
\scalebox{0.35}{{\includegraphics{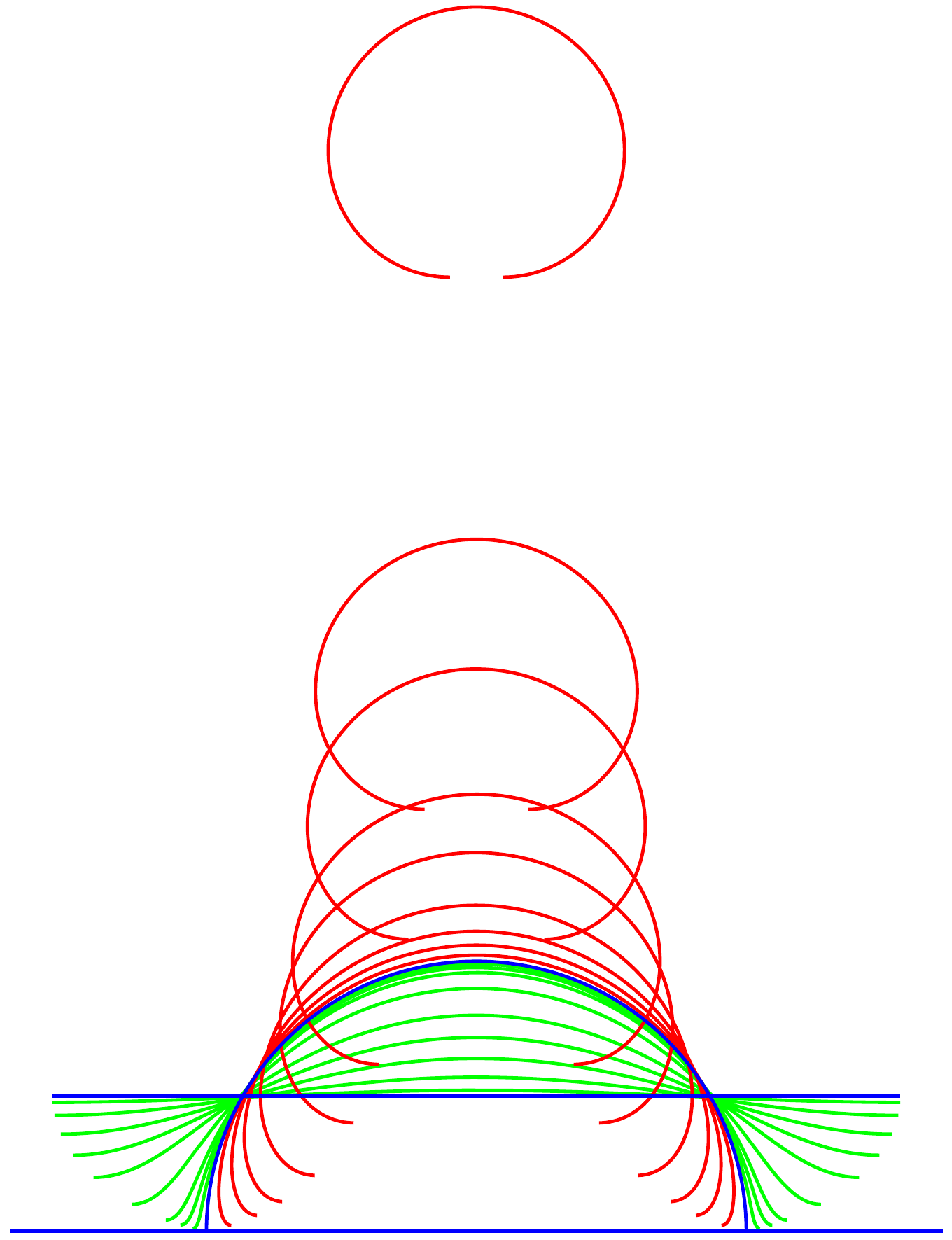}}}\hspace{.1cm}\scalebox{0.35}
{{\includegraphics{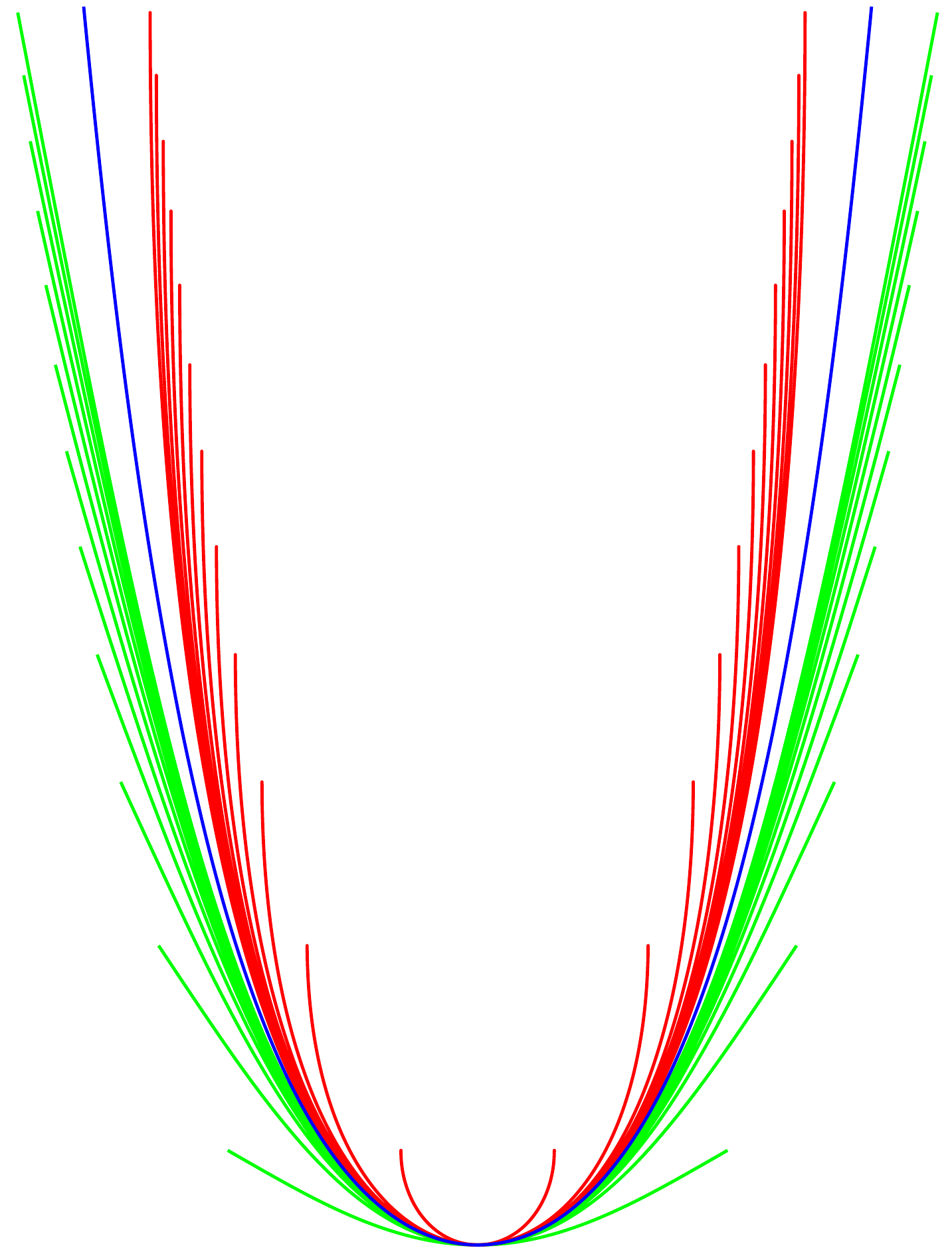}}}
\end{center}
\caption{\label{familia2} Left: Location of the roulettes in Figure \ref{familia1} with respect to the abscissa axis.
Right: a family of roulettes of type $B_1$ (green) and $C_1$ (red) separated by a catenary (blue).}
\end{figure}
 
In Figure \ref{familia2} (left) we show the same family of roulettes shown in Figure \ref{familia1}, but to the height, 
with respect to the abscissa axis, at which it has been generated by the focus of the corresponding conic. Observe that 
the limiting circumference of radius $a$ would have the center at height $b=\infty$. In Figure \ref{familia2}(right) we 
display the asymptotic role of the catenary separating two families of curves $B_1$ and $C_1$. In this case the corresponding 
conics don't have the parameter $a$ constant, and both families tends to the catenary as the parameter $a$ grows.

\section{The Delaunay Surfaces} In this section we study Delaunay surfaces and derive analytical expressions for their most 
important differential geometric properties. The Delaunay surfaces are surfaces of revolution and therefore the key to their properties 
lie in their meridians, which here are the roulettes of the foci of the conics discussed in the previous section.
The Delaunay surfaces are thus the surfaces of revolution generated by the curves  $A(t)$, $B_1(t)$, $B_2(t)$, $C_1(t)$ 
and $C_2(t)$. We next describe the parametrization, the coefficients of the first and second fundamental forms of each one of 
these surfaces and their curvatures. The entire differential structure has a ramarkably transparent dependence on the 
parameters that characterize each conic.  

\begin{figure}[htb]
\begin{center}
\scalebox{0.18}{{\includegraphics{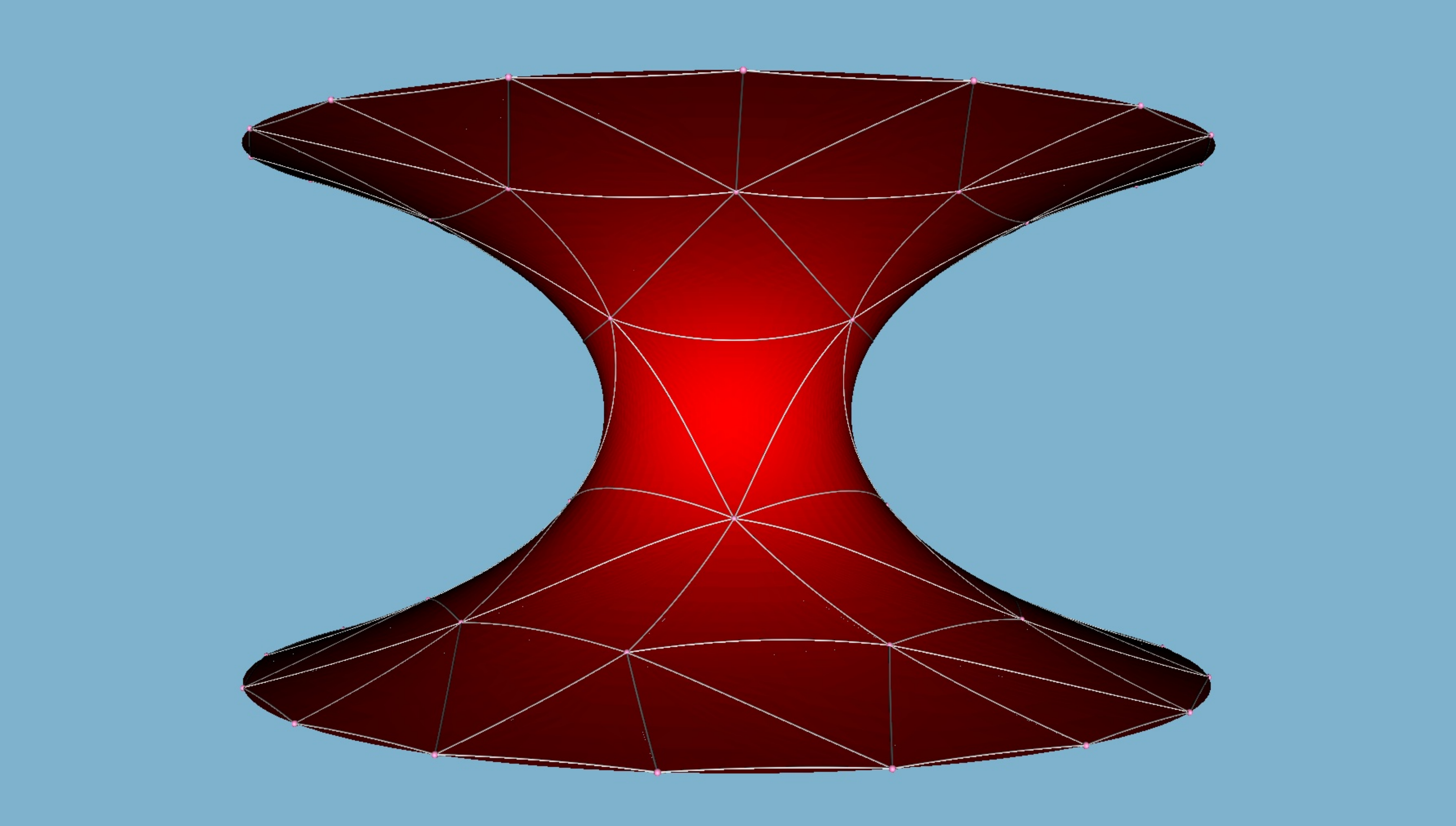}}}
\end{center}
\caption{\label{catenoid} The catenoid.}
\end{figure}

{\sc Catenoid:} $\pmb x(t,v)=\Big(f_c(t)\cos(v),f_c(t)\sin(v),g_c(t)\Big),$ {--} see Figure \ref{catenoid}. We have
$$\begin{array}{rl}
 {\pmb x}_t=&\hspace{-.2cm}\left(b\sinh(t)\cos(v),
b\sinh(t)  \sin(v),b  \right), \\[1ex]
 {\pmb x}_v=&\hspace{-.2cm}\left(-b\cosh(t)\sin(v), b\cosh(t)\cos(v),0\right).  \end{array}$$
The unit normal vector at $(t,v)$ is given by
$${\pmb n}_c(t,v)=\left(-\frac{\cos(v)}{\cosh(t)},-\frac{\sin(v)}{\cosh(t)},\tanh(t)\right).$$
The non--vanishing coeficients of the first and the second fundamental form, and the principal curvatures, the mean curvature 
and the Gaussian curvature are,  
$$\begin{array}{rlrlrlrl}
E=&\hspace{-.25cm}\displaystyle b^2\cosh^2(t),&  G=&\hspace{-.25cm}\displaystyle b^2\cosh^2(t),&
L=&\hspace{-.25cm}\displaystyle -b,&N=&\hspace{-.25cm}\displaystyle b,\\[1ex]
k_1=&\hspace{-.25cm}\displaystyle\frac{-1}{b\cosh^2(t)},&
k_2=&\hspace{-.25cm}\displaystyle \frac{1}{b\cosh^2(t)},&
 H= &\hspace{-.25cm}\displaystyle  0,& K= &\hspace{-.25cm}\displaystyle\frac{-1}{b^2\cosh^4(t)}.
 \end{array} $$
The geodesic curvature of a parallel is
$$k_g=\frac{\sinh(t)}{b\cosh^2(t)},$$
and the total curvature of the catenoid is
$$\int_{\pmb x} K d\sigma=
-(v_2-v_1) \tanh(t)\Big|_{t_1}^{t_2}.$$

\begin{figure}[htb]
\begin{center}
\scalebox{0.25}{{\includegraphics{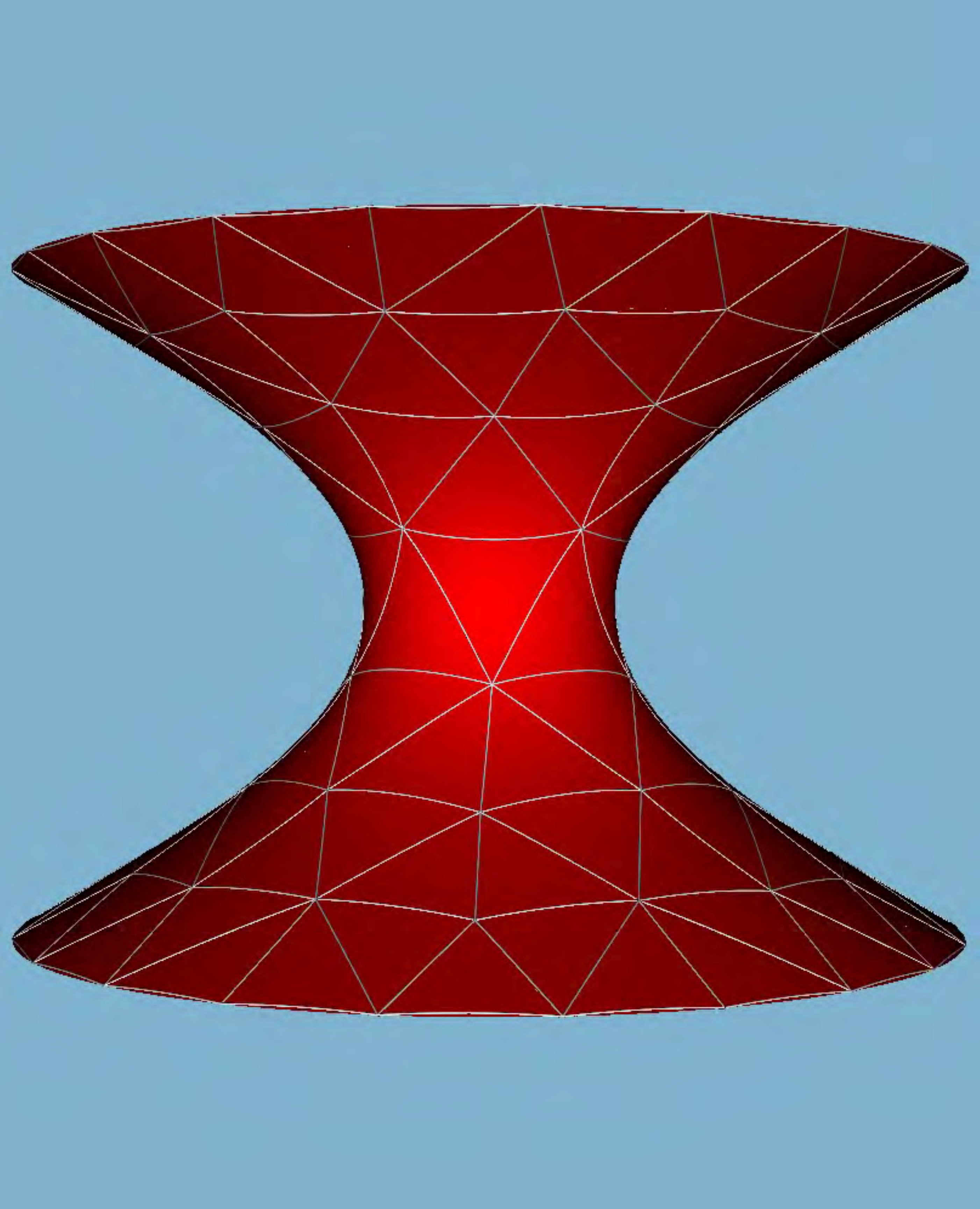}}}
\scalebox{0.25}{{\includegraphics{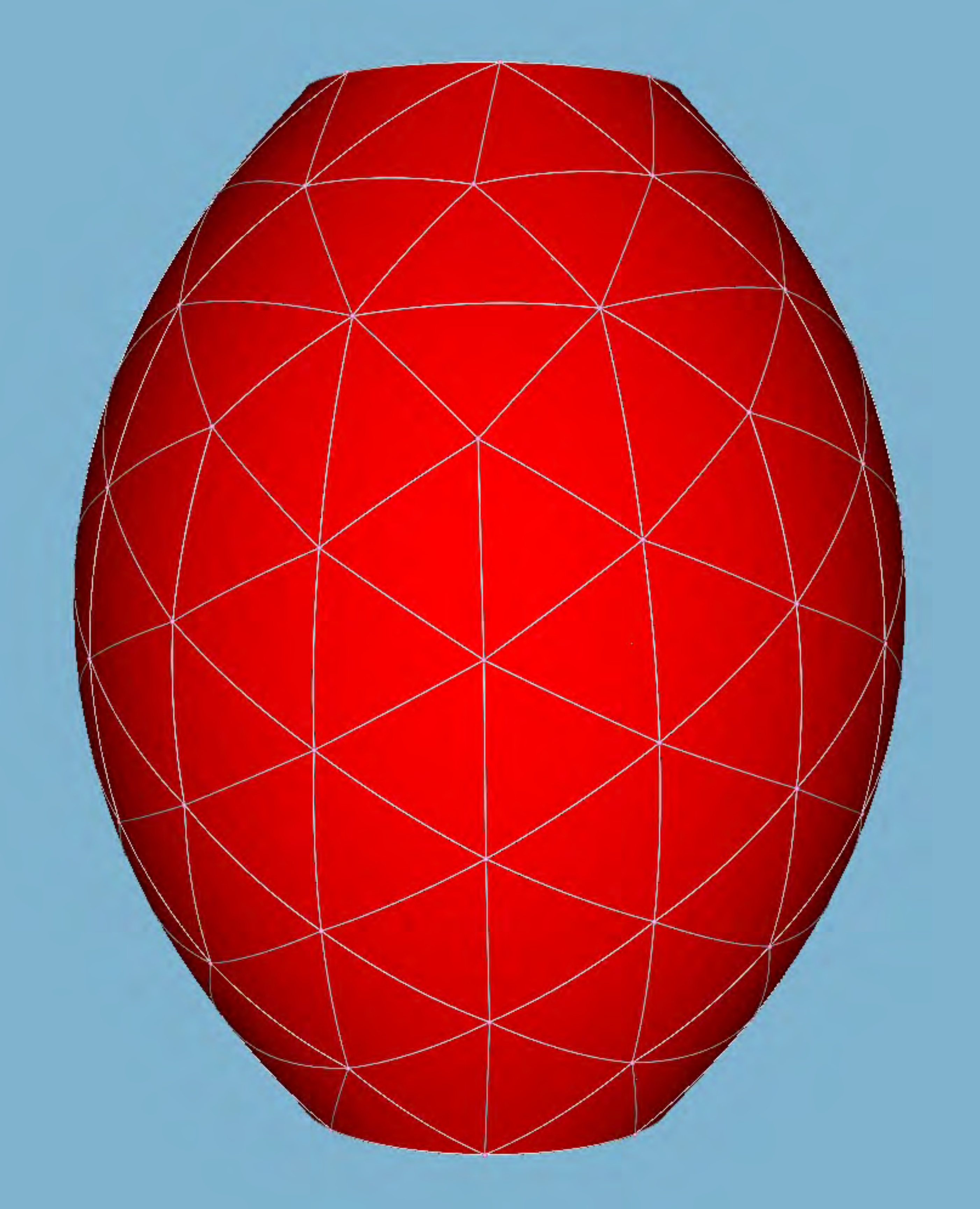}}}
\end{center}
\caption{\label{unduloids} Unduloids generated by the revolution of $B_1$ (left), and $B_2$ (right).}
\end{figure}

{\sc Unduloid:} ${\pmb y}(t,v)=\Big(f_u^1(t)\cos(v),f_u^1(t)\sin(v),g_u^1(t)\Big),$ {--} see Figure \ref{unduloids} (left).
Considering
$$h_u(t)=\frac{ab}{\sqrt{a^2-c^2\cos^2(t)}\left(a+c\cos(t)\right)}$$
it follows that
$$\begin{array}{rl}
 {\pmb y}_t=&\hspace{-.2cm}\left(c\sin(t) h_u(t)\cos(v),
c\sin(t) h_u(t) \sin(v),b h_u(t)\right), \\[1ex]
 {\pmb y}_v=&\hspace{-.2cm}\left(-f_u^1(t)\sin(v), f_u^1(t)\cos(v),0\right).  \end{array}$$
The unit normal vector at $(t,v)$ is given by
$${\pmb n}_u(t,v)=\left(\frac{-b\cos(v)}{\sqrt{a^2-c^2\cos^2(t)}},
\frac{-b\sin(v)}{\sqrt{a^2-c^2\cos^2(t)}},
\frac{c\sin(t)}{\sqrt{a^2-c^2\cos^2(t)}}\right).$$

$$\begin{array}{rlrl }
E=&\hspace{-.25cm}\displaystyle \frac{a^2b^2}{\left(a+c\cos(t)\right)^2},&\hspace{1cm} G=&\hspace{-.25cm}\displaystyle \frac{b^2\left(a-c\cos(t)\right)}{\left( a+c\cos(t)\right)},\\[3ex]
L=&\hspace{-.25cm}\displaystyle \displaystyle\frac{-ab^2c\cos(t)}{\left( a^2-c^2\cos^2(t)   \right)\left(a+c\cos(t)\right)},&N=&\hspace{-.25cm}\displaystyle \frac{b^2}{a+c\cos(t)},\\[3ex]
k_1=&\hspace{-.25cm}\displaystyle\frac{-c\cos(t)}{a(a-c\cos(t))},&
k_2=&\hspace{-.25cm}\displaystyle \frac{1}{a-c\cos(t)},\\[3ex]
 H= &\hspace{-.25cm}\displaystyle  \frac{1}{2a},& K= &\hspace{-.25cm}\displaystyle\frac{-c\cos(t)}{\displaystyle a\left(a-c\cos(t)\right)^2}.
 \end{array} $$
The geodesic curvature of a parallel is
$$k_g=\frac{-c\sin(t)}{b\big(a-c\cos(t)\big)},$$
and the total curvature is
$$\int_{\pmb y} K d\sigma= \left.
-c(v_2-v_1) \frac{\sin(t)}{\sqrt{a^2-c^2\cos^2(t)}} \right|_{t_1}^{t_2}.$$

Clearly one could do the same calculation for the surface generated by $B_2$ (see Fig \ref{unduloids}(right)), but it is enough  
consider the proper domain of the parameter $t$ to belong to one or another part of the unduloid. In fact these 
surfaces are periodic with period $2\pi$ with respect to the parameter $t$. 

We must, however, consider both parts $C_1$ and $C_2$ in the construction of the nodoids, because each roulette 
has its domain in $\RR$.

\begin{figure}[htb]
\begin{center}
\scalebox{0.13}{{\includegraphics{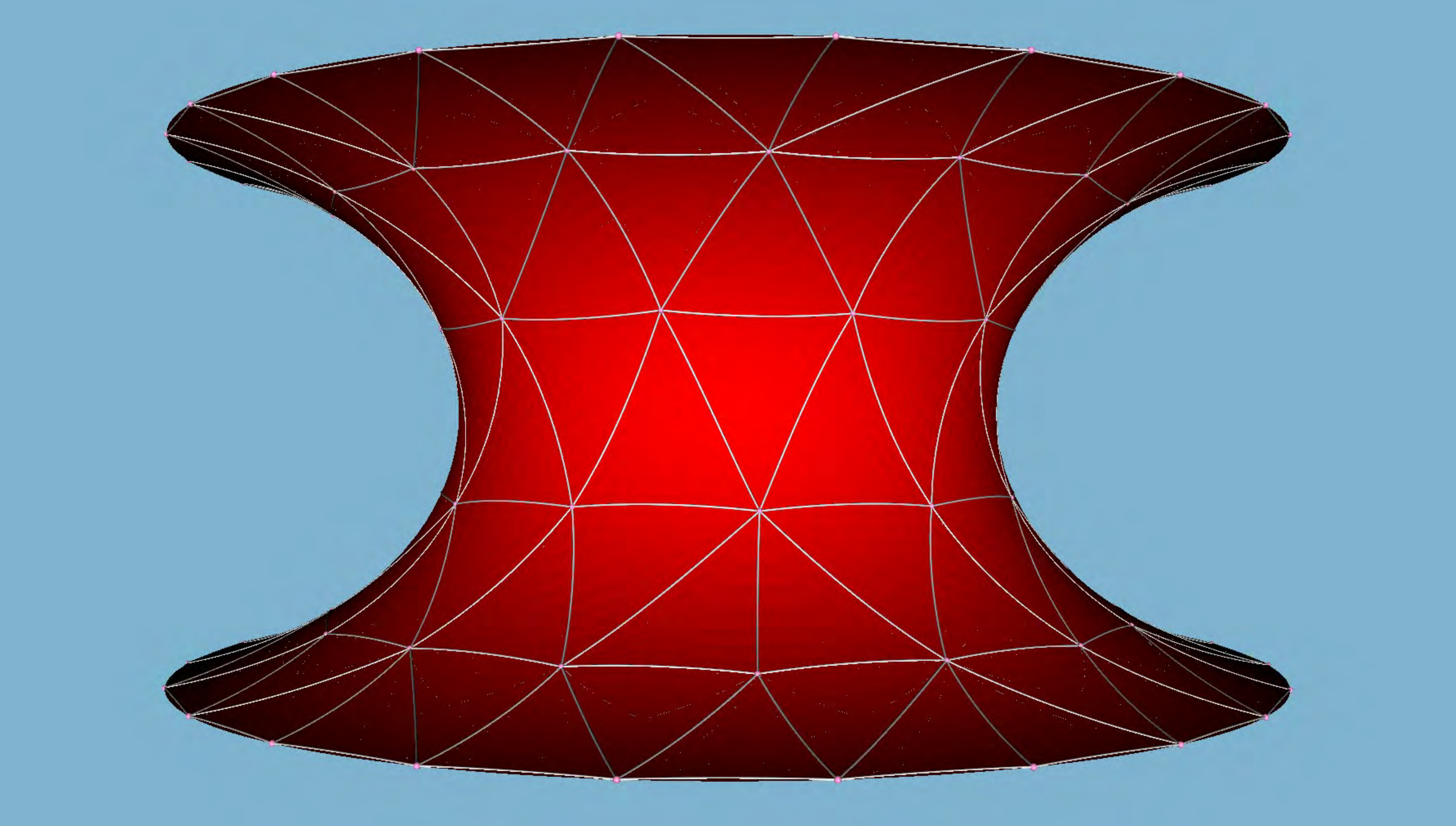}}}
\scalebox{0.13}{{\includegraphics{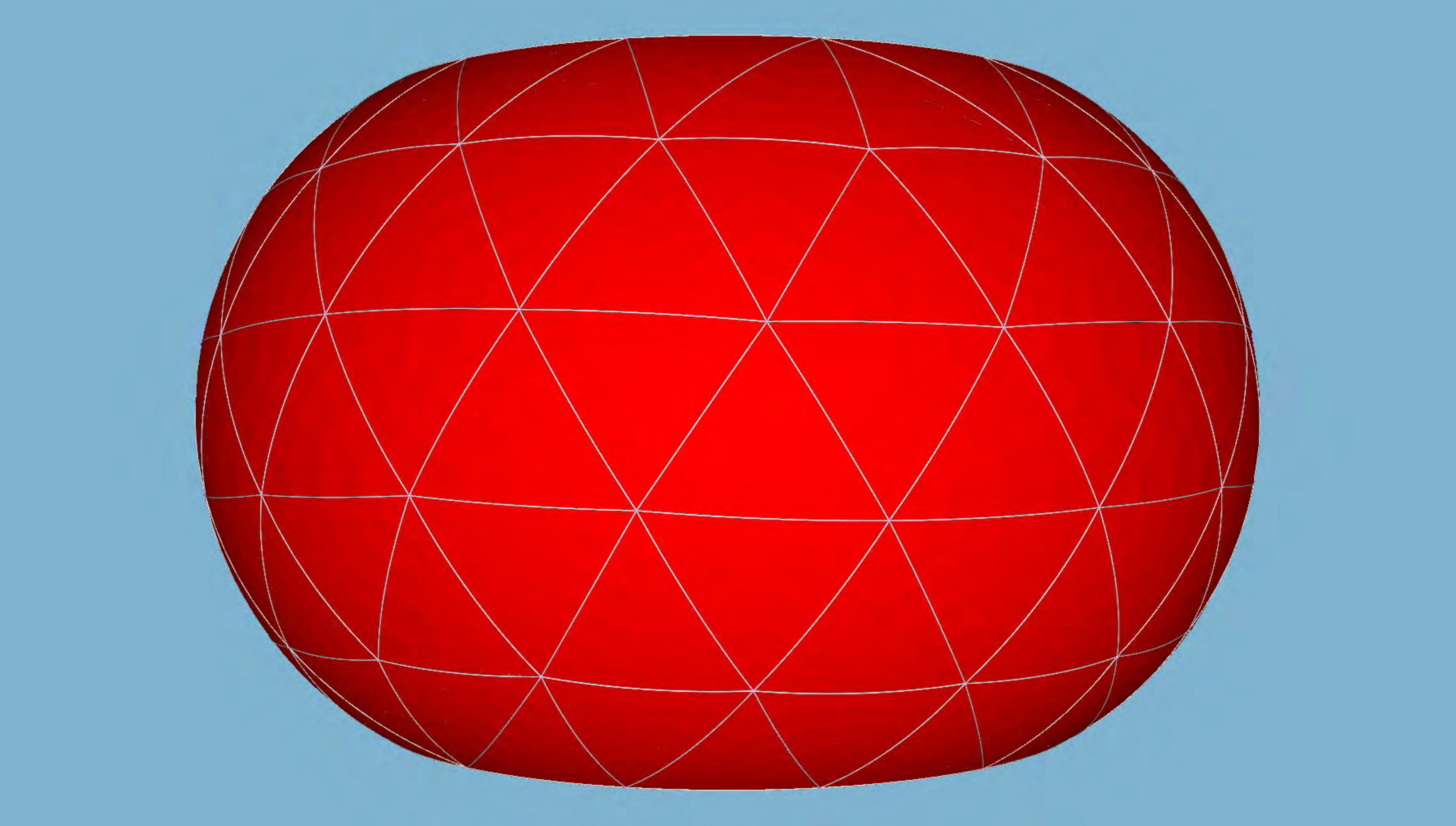}}}
\end{center}
\caption{\label{nodoids} Nodoids generated by the revolution of $C_1$ (left) and $C_2$ (right).}
\end{figure}

{\sc Nodoid1:} $\pmb z_1(t,v)=\Big(f_n^1(t)\cos(v),f_n^1(t)\sin(v),g_n^1(t)\Big),$ {--} see Figure \ref{nodoids}(left). 
We consider
$$h_1(t)=\frac{ab}{\sqrt{c^2\cosh^2(t)-a^2}\left(c\cosh(t)+a\right)}$$
from which it follows that
$$\begin{array}{rl}
 ({\pmb z}_1)_t=&\hspace{-.2cm}\Big(c\sinh(t) h_1(t)\cos(v),
c\sinh(t) h_1(t) \sin(v),b h_1(t)\Big), \\[1ex]
 ({\pmb z}_1)_v=&\hspace{-.2cm}\Big(-f_n^1(t)\sin(v), f_n^1(t)\cos(v),0\Big).  \end{array}$$
The unit normal vector at $(t,v)$ is given by
$${\pmb n}_1(t,v)=\left(\frac{-b\cos(v)}{\sqrt{c^2\cosh^2(t)-a^2}},
\frac{-b\sin(v)}{\sqrt{c^2\cosh^2(t)-a^2}},
\frac{c\sinh(t)}{\sqrt{c^2\cosh^2(t)-a^2}}\right).$$

$$\begin{array}{rlrl }
E=&\hspace{-.25cm}\displaystyle \frac{a^2b^2}{\left(c\cosh(t)+a\right)^2},&\hspace{1cm} G=&\hspace{-.25cm}\displaystyle \frac{b^2\left(c\cosh(t)-a\right)}{ c\cosh(t)+a},\\[3ex]
L=&\hspace{-.25cm}\displaystyle \frac{-ab^2c\cosh(t)}{\left( c^2\cosh^2(t)-a^2   \right)\left(c\cosh(t)+a\right)},&N=&\hspace{-.25cm}\displaystyle \frac{b^2}{c\cosh(t)+a},\\[3ex]
k_1=&\hspace{-.25cm}\displaystyle\frac{-c\cosh(t)}{\displaystyle a(c\cosh(t)-a)},&
k_2=&\hspace{-.25cm}\displaystyle \frac{1}{\displaystyle c\cosh(t)-a},\\[3ex]
 H= &\hspace{-.25cm}\displaystyle  \frac{-1}{2a},& K= &\hspace{-.25cm}\displaystyle\frac{-c\cosh(t)}{\displaystyle a\left(c\cosh(t)-a\right)^2}.
 \end{array} $$
The geodesic curvature of a parallel is
$$k_g=\frac{c\sinh(t)}{b\big(c\cosh(t)-a\big)},$$
and the total curvature is
$$\int_{\pmb z_1} K d\sigma= -c(v_2-v_1) \left.\frac{\sinh(t)}{\sqrt{c^2\cosh^2(t)-a^2}} \right|_{t_1}^{t_2}.$$

\begin{figure}[htb]
\begin{center}
\scalebox{0.12}{{\includegraphics{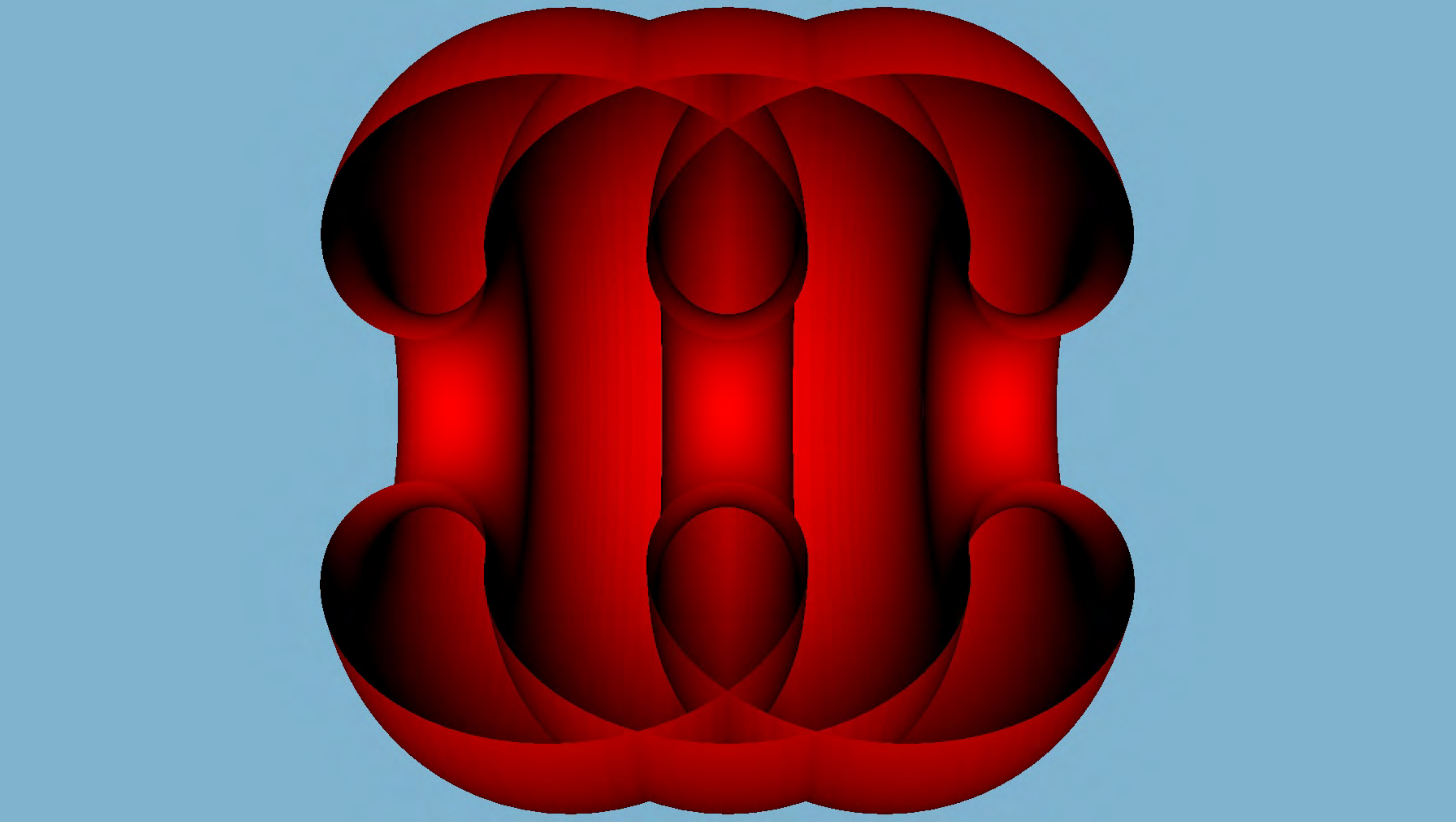}}}
\scalebox{0.12}{{\includegraphics{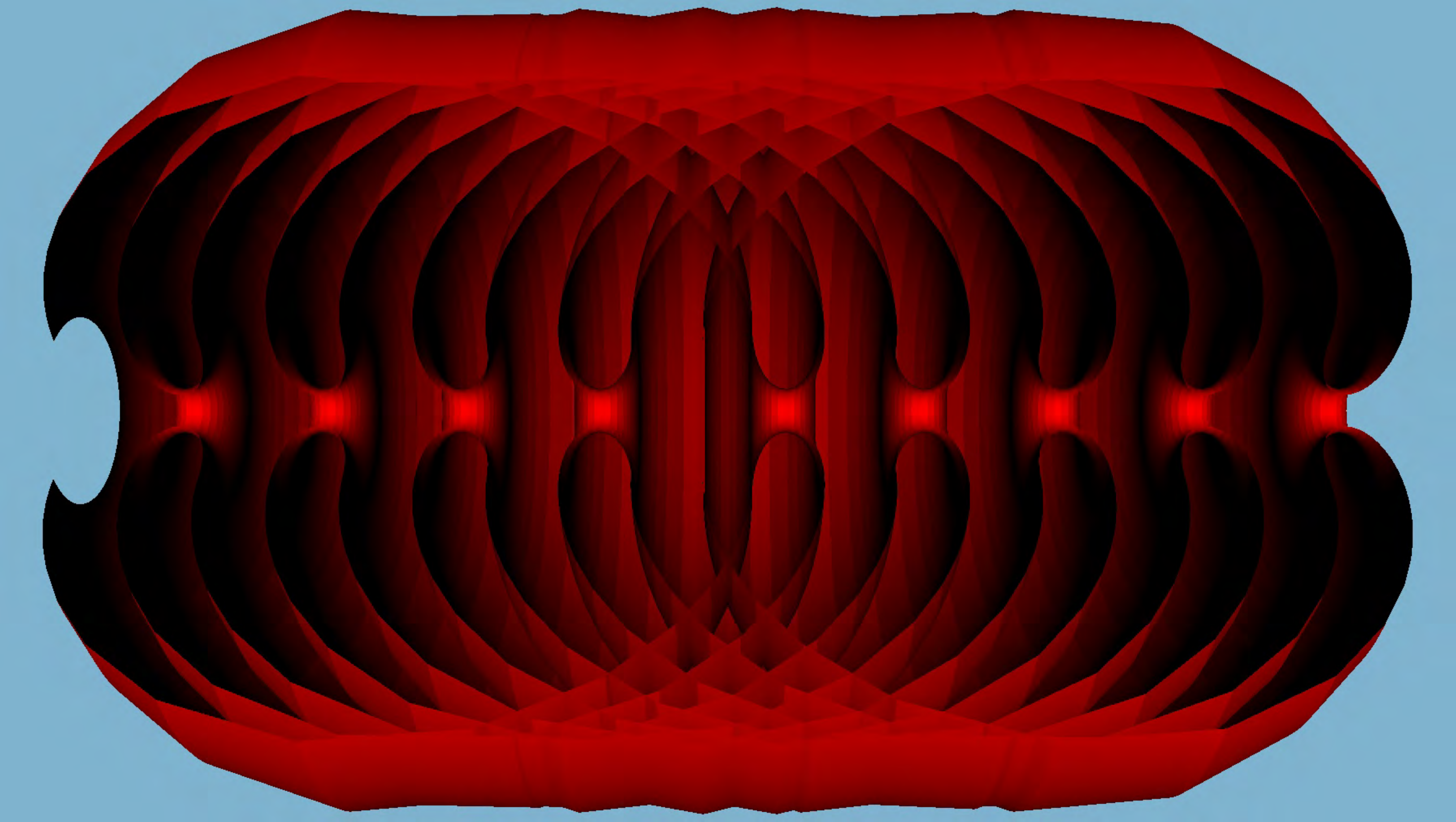}}}
\end{center}
\caption{\label{nodoidsurfaces} Left: Closed nodoid (compact without boundary), connected sum of 4 tori. 
Right: Connected compact nodoid with boundary.}
\end{figure}

{\sc Nodoid2:} $\pmb z_2(t,v)=\Big(f_2(t)\cos(v),f_2(t)\sin(v),g_2(t)\Big),$ {--} see Figure \ref{nodoids}(right). 
We consider
$$h_2(t)=\frac{-ab}{\sqrt{c^2\cosh^2(t)-a^2}\left(c\cosh(t)-a\right)}$$
from which it follows that
$$\begin{array}{rl}
 ({\pmb z}_2)_t=&\hspace{-.2cm}\Big(c\sinh(t) h_2(t)\cos(v),
c\sinh(t) h_2(t) \sin(v),b h_2(t)\Big), \\[1ex]
 ({\pmb z}_2)_v=&\hspace{-.2cm}\Big(-f_2(t)\sin(v), f_2(t)\cos(v),0\Big).  \end{array}$$
The unit normal vector at $(t,v)$ is given by
$${\pmb n}_2(t,v)=\left(\frac{b\cos(v)}{\sqrt{c^2\cosh^2(t)-a^2}},
\frac{b\sin(v)}{\sqrt{c^2\cosh^2(t)-a^2}},
-\frac{c\sinh(t)}{\sqrt{c^2\cosh^2(t)-a^2}}\right).$$

$$\begin{array}{rlrl }
E=&\hspace{-.25cm}\displaystyle \frac{a^2b^2}{\left(c\cosh(t)-a\right)^2},&\hspace{1cm} G=&\hspace{-.25cm}\displaystyle \frac{b^2\left(c\cosh(t)+a\right)}{c\cosh(t)-a},\\[3ex]
L=&\hspace{-.25cm}\displaystyle \frac{-ab^2c\cosh(t)}{\left( c^2\cosh^2(t)-a^2   \right)\left(c\cosh(t)-a\right)},&N=&\hspace{-.25cm}\displaystyle \displaystyle\frac{-b^2}{c\cosh(t)-a},\\[3ex]
k_1=&\hspace{-.25cm}\displaystyle\frac{-c\cosh(t)}{\displaystyle a(c\cosh(t)+a)},&
k_2=&\hspace{-.25cm}\displaystyle\frac{-1}{\displaystyle c\cosh(t)+a},\\[3ex]
 H= &\hspace{-.25cm}\displaystyle  \frac{-1}{2a},& K= &\hspace{-.25cm}\displaystyle\frac{c\cosh(t)}{\displaystyle a\left(c\cosh(t)+a\right)^2}.
 \end{array} $$
The geodesic curvature of a parallel is
$$k_g=\frac{-c\sinh(t)}{b\big(c\cosh(t)+a\big)}.$$
The total curvature is
$$\int_{\pmb z_2} K d\sigma= c(v_2-v_1)\left.\frac{\sinh(t)}{\sqrt{c^2\cosh^2(t)-a^2}} \right|_{t_1}^{t_2}.$$

\begin{figure}[htb]
\begin{center}
\scalebox{0.6}{{\includegraphics{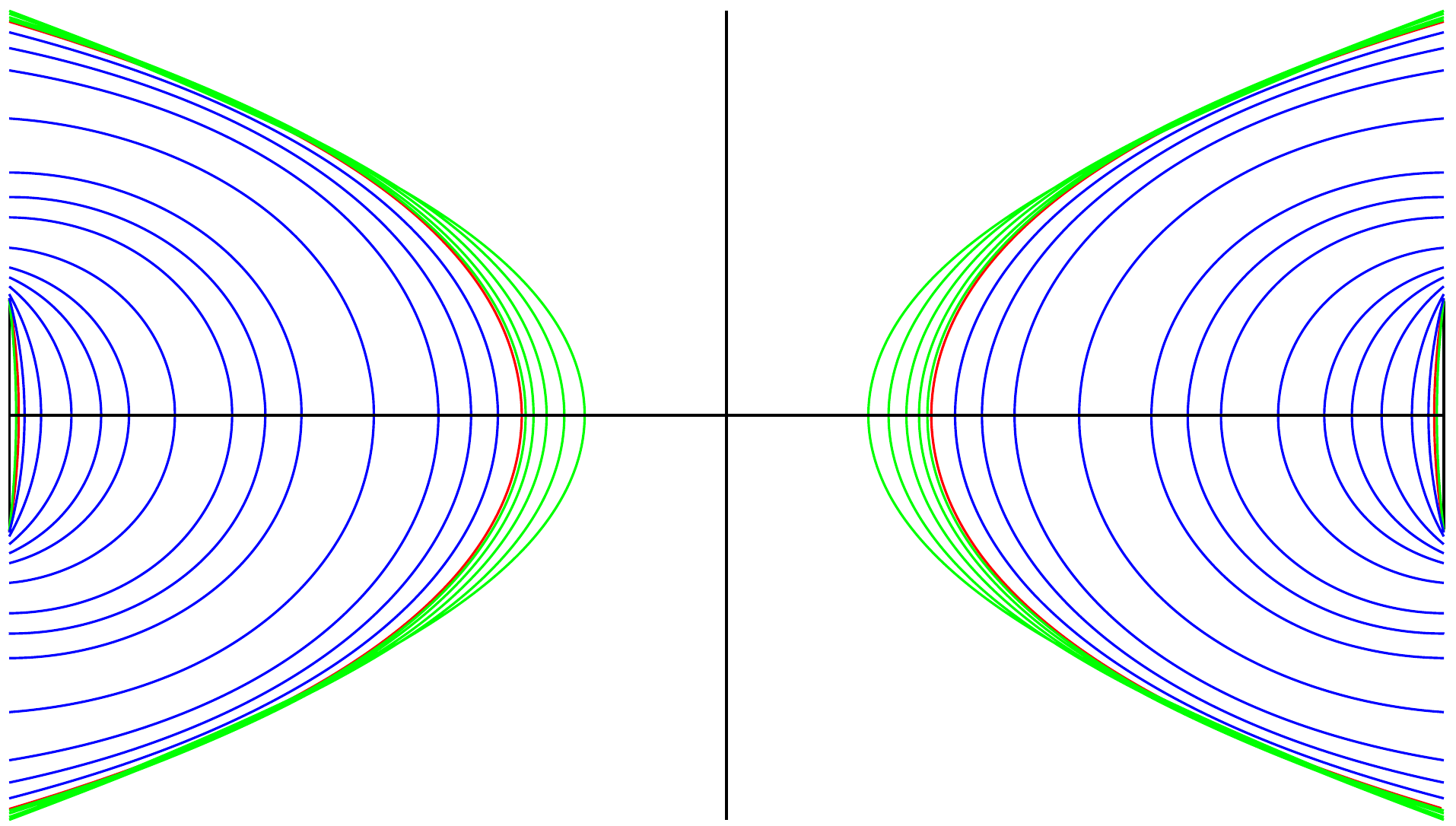}}}
\end{center}
\caption{\label{volumen} Planar sections containing the revolution axis of a family of Delaunay surfaces with constant volume. 
From the outside-in: cylinder (black), unduloids (green), catenoid (red), nodoids (blue), catenoid (red) and unduloids (green).}
\end{figure}

Fig. \ref{nodoidsurfaces} illustrates the versatility of the parameterization of nodoids adopted here. 
We avoid the periodic extension of nodoids with both positive and negative Gaussian curvature, 
maintain the orientability. preserve the $C^\infty$ class and find new surfaces both closed and 
with boundary.

Delaunay surfaces are characterized by minimizing the area with fixed boundaries and constant volume. Here we illustrate the 
versatility of our formulation by finding a surface which satisfies these conditions. Consider a plane curve $(f(t), g(t))$.
The volume enclosed by its surface of revolution is given by 
$\pi\int_{t_0}^{t_1} f^2(t) g'(t)dt$ and the end points by $\left( f(t_0), g(t_0)\right)=(f_0, g_0)$ and 
$\left( f(t_1), g(t_1)\right)=(f_1, g_1)$. 

Consider for example the symmetric nodoid generated by a roulette $C_1$ such that the enclosed 
volume is 1 and the radius is 1 at the end points $\pm t_0$. The equations to solve for $a$ and $b$ are then 

$$\begin{array}{rlrl }
\displaystyle\pi a b^4\int_{-t_0}^{t_0}\frac{\left( c\cosh(t)-a\right)}{\left( c\cosh(t)+a\right)^2\sqrt{c^2\cosh^2(t)-a^2}}dt&=1,\\[5ex]
\displaystyle\frac{b\left(c\cosh(t_0)-a\right)}{\sqrt{c^2\cosh^2(t_0)-a^2}}&=1.
 \end{array} $$ 
In Fig \ref{volumen} we present planar sections that contain the common axis of revolution of several Delaunay surfaces with enclosed 
volume equal to 1 and radius 1 at the boundaries. Note that the solutions depend on $t_0$. 

Finally the paremetrizations developed here have proven extremely useful in analytic and computational 
explorations of the structure of crystalline particle arrays on Delaunay surfaces as realized experimentally 
in capillary bridges \cite{BBMY13}.

\noindent{\bf Acknowledgments.}  This work has been partly supported by the Spanish Research Council (Comisi\'on
Interministerial de Ciencia y Tecnolog\'\i a,) under project MTM2010-19660. The work of M.J.B. was also supported by the 
National Science Foundation through Grant No. DMR-1004789 and by funds from the Soft Matter Program of Syracuse University.

\end{document}